\input amstex
\input colordvi
\documentstyle{amsppt}  

\pagewidth{12.5cm}\pageheight{19cm}\magnification\magstep1
\topmatter
\title Coordinate rings and birational charts \endtitle
\author Sergey Fomin and George Lusztig\endauthor
\address Department of Mathematics, M.I.T., Cambridge, MA 02139, USA \endaddress
\address Department of Mathematics, University of Michigan, Ann Arbor, MI 48109, USA \endaddress
\abstract Let $G$ be a semisimple simply connected complex algebraic group.
Let $U$ be the unipotent radical of a Borel subgroup in~$G$. 
We describe the coordinate rings of $U$ (resp., $G/U$, $G$) 
in terms of two (resp., four, eight) birational charts introduced in [L94, L19] in connection with the study of total positivity.
\endabstract
\thanks{This work was supported by NSF grants DMS-1664722 (S.~F.) and 
DMS-1855773 (G.~L.) and by a Simons Fellowship (S.~F.).}
\endthanks
\endtopmatter   

\document

\define\dz{\dot z}

\define\dw{\dot w}

\define\ds{\dot s}
\define\dy{\dot y}

\define\frl{\forall}

\define\si{\sim}

\define\sqc{\sqcup}

\define\part{\partial}
\define\emp{\varnothing}

\define\ra{\rangle}
\define\n{\notin}
\define\iy{\infty}
\define\m{\mapsto}
\define\do{\dots}
\define\la{\langle}

\define\sub{\subset}    

\define\T{\times}
\define\ti{\tilde}
\define\nl{\newline}
\redefine\i{^{-1}}

\define\ov{\overline}

\define\Hom{\text{\rm Hom}}

\define\a{\alpha}

\redefine\c{\chi}
\define\g{\gamma}
\redefine\d{\delta}
\define\e{\varepsilon}

\define\io{\iota}
\redefine\o{\omega}

\define\ph{\phi}

\redefine\l{\lambda}

\redefine\D{\Delta}

\define\ii{\bold i}
\define\jj{\bold j}

\define\CC{\Bbb C}

\define\FF{\Bbb F}

\define\RR{\Bbb R}

\define\ZZ{\Bbb Z}

\define\cb{\Cal B}

\define\car{\Cal R}

\define\cu{\Cal U}
\define\cv{\Cal V}

\define\cz{\Cal Z}
\define\cx{\Cal X}
\define\cy{\Cal Y}

\define\sha{\sharp}

\define\che{\check}
\define\cha{\che{\a}}

\head Introduction\endhead
Let $G$ be a simply connected, almost simple algebraic group over~$\CC$. 
Fix a maximal torus $T$ of $G$ and a pair $B^+,B^-$ 
of opposite Borel subgroups containing~$T$, with unipotent radicals $U^+,U^-$. 
Let $\nu=\dim(U^+)$ and $r=\dim(T)$.
For an irreducible quasi-affine variety $X$ over $\CC$, we denote by $O(X)$ the 
algebra of regular functions $X@>>>\CC$, and let $[O(X)]$ be the quotient field of~$O(X)$.

In this paper, we show (see Theorems 0.3, 4.2 and 5.2) that the algebra
$O(U^+)$ (resp., $O(G/U^-)$ and $O(G)$) can be completely described
in terms of two (resp., four and eight) birational charts
$\CC^\nu@>>>U^+$ (resp., $\CC^\nu\T(\CC^*)^r@>>>G/U^-$ and
$\CC^\nu\T(\CC^*)^r\T\CC^\nu@>>>G$)
which were introduced in \cite{L94},\cite{L19} in connection with the study of total positivity.

Theorem 0.3 provides a proof of a conjecture in \cite{L19, 6.1(a)}.
Theorem~4.2 (resp., Theorem~5.2) establishes a weak form of a
conjecture in \cite{L19,~6.3(a)} (resp., \cite{L19, 6.2(a)}) in which
only two birational charts, instead of four (resp., eight), were used.
The proof of Theorem~0.3 given in Section 3 relies on the results in \cite{BZ97} and \cite{FZ99}
that describe the inverse of the charts for $U^+$ in terms
of ``generalized minors.'' 
Theorems 4.2 and 5.2 are proved in Sections 4 and~5, respectively, using reduction
to the case of~$U^+$.   In particular, our proof of Theorem 5.2 does not
use the more complete results on generalized minors in~\cite{FZ99}. 
(The latter technique would have allowed to decrease the number of charts from eight to two, 
but then the two charts used would not be canonical, unlike the eight that we consider~here.) 

\pagebreak

In order to state our main result (Theorem~0.3), we will need to introduce some notation. 

Let $U_i^+$ ($i\in I$) be the simple root subgroups of $U^+$,
and let $U_i^-$ ($i\in I)$ be the corresponding root
subgroups of $U^-$; here $I$ is a finite indexing set.
We assume that for any $i\in I$ we are given isomorphisms of
algebraic groups $x_i:\CC@>\si>>U_i^+$ and $y_i:\CC@>\si>>U_i^-$ such that
 $(T,B^+,B^-,x_i,y_i;i\in I)$ is a pinning for $G$.
 
\definition{Definition~0.1}
Let $I^*$ be the set of all pairs $(i,j)\in I\T I$ such that
any element in $U_i$ commutes with any element in $U_j$.
There is a unique (up to a labeling convention) partition $I=I_0\sqc I_1$ into two disjoint subsets such that
$I_0\T I_0\sub I^*$  and $I_1\T I_1\sub I^*$. 
Let $r_0=\sha(I_0)$ and $r_1=\sha(I_1)$ be the cardinalities of $I_0$ and~$I_1$.

It is known that $h=2\nu/r$ is an integer (the Coxeter number). 

For $\e\in\ZZ$, we define $[\e]\in\{0,1\}$ by $\e\equiv[\e]\bmod2$.
With this notation, we have 
$$
\nu=\underbrace{r_{[\e]}+r_{[\e+1]}+\do+r_{[\e+h-1]}}_{\text{$h$ terms}}.$$
(If $h$ is even, this follows from $r_0+r_1=r$; if $h$ is odd, we use 
that $r_0=r_1=r/2$.) 

For $\e\in\{0,1\}$, let us fix the ordering of the elements of $I_\e$: 
$$I_\e=\{i_1^\e,i_2^\e,\do,i_{r_\e}^\e\}.$$
We then define the sequence $\jj^\e\in I^\nu$  (a distinguished reduced expression) by
$$
\aligned
\jj^\e&=(j^\e_1,j^\e_2,\do,j^\e_\nu) \\
&=(i_1^{[\e]},i_2^{[\e]},\do,i_{r_{[\e]}}^{[\e]},
i_1^{[\e+1]},i_2^{[\e+1]},\do,i_{r_{[\e+1]}}^{[\e+1]}, \\
&\qquad 
i_1^{[\e+2]},i_2^{[\e+2]},\do,i_{r_{[\e+2]}}^{[\e+2]},
\do,
i_1^{[\e+h-1]},i_2^{[\e+h-1]},\do,i_{r_{[\e+h-1]}}^{[\e+h-1]}).\endaligned \tag 0.1.1
$$
(The upper indices are not exponents.)
Thus, the first $r_{[\e]}$ terms of $\jj^\e$ are the elements of
$I_{[\e]}$ in their order, the next $r_{[\e+1]}$
terms are the elements of $I_{[\e+1]}$ in their order,
and these patterns keep alternating until we accumulate $\nu$ entries.

For a sequence of indices $\ii=(i_1,i_2,\do,i_n)\in I^n$ of length $n\ge0$, we define
the map $f_\ii:\CC^n@>>>U^+$ by
$$f_\ii(a_1,a_2,\do,a_n)=x_{i_1}(a_1)x_{i_2}(a_2)\do x_{i_n}(a_n).
\tag 0.1.2$$
\enddefinition

In particular, one can choose $\ii=\jj^\e$ for $\e\in\{0,1\}$, as in (0.1.1) above. 
The following fact is well known:

\proclaim{Proposition 0.2} 
The maps $f_{\jj^0},f_{\jj^1}$ are birational isomorphisms from $\CC^\nu$ to
$U^+$. 
\endproclaim

Proposition 0.2 can be deduced from the proof of \cite{L94, 2.7} using (1.3.1) below;
it can also be deduced from \cite{BZ97}. See also 3.12(d).

By Proposition 0.2, each map $f_{\jj^\e}$ ($\e\in\{0,1\}$) induces an isomorphism of fields 
$f_{\jj^\e}^*:[O(U^+)]@>\si>>[O(\CC^\nu)]$.

\proclaim{Theorem 0.3} An element $\ph\in[O(U^+)]$ belongs to $O(U^+)$
if and only if the rational function $f_{\jj^\e}^*(\ph)\in[O(\CC^\nu)]$ belongs to $O(\CC^\nu)$
for $\e=0$ and for $\e=1$.
\endproclaim
The proof of Theorem 0.3 is given in Section~3.

The instances of Theorem~0.3 for $G$ of types $A_2$ and $A_3$ have been verified
in \cite{L19, Section~6.1}. 
In the rest of this section, we work out the latter case in detail.  

\definition{Example 0.4}
Let $G=\operatorname{SL}_4(\CC)$, with $T$, $B^+$, and $B^-$ 
its subgroups of diagonal, upper-triangular, and low-triangular matrices,
respectively. 
Then 
$$
\align
U^+&=\left\{ 
\left.
u=\left[
\matrix
1 & u_{12} & u_{13} & u_{14} \\
0 & 1 & u_{23} & u_{24} \\
0 & 0 & 1 & u_{34} \\
0 & 0 & 0 & 1
\endmatrix
\right]
\right| u_{12}, u_{13}, u_{14}, u_{23}, u_{24}, u_{34}\in\CC
\right\}, \\
O(U^+)&=\CC[u_{12}, u_{13}, u_{14}, u_{23}, u_{24}, u_{34}], \\
r&=3, \\
I&=\{1,2,3\}, \\
x_1(a)&\!=\!\left[
\matrix
1 & a & 0 & 0 \\
0 & 1 & 0 & 0 \\
0 & 0 & 1 & 0 \\
0 & 0 & 0 & 1
\endmatrix
\right]\!,
\ \ x_2(a)\!=\!\left[
\matrix
1 & 0 & 0 & 0 \\
0 & 1 & a & 0 \\
0 & 0 & 1 & 0 \\
0 & 0 & 0 & 1
\endmatrix
\right]\!,
\ \ x_3(a)\!=\!\left[
\matrix
1 & 0 & 0 & 0 \\
0 & 1 & 0 & 0 \\
0 & 0 & 1 & a \\
0 & 0 & 0 & 1
\endmatrix
\right]\!,
\\
\nu&=6, \\
h&=4. 
\endalign
$$
We set
$I_0=\{2\}$ and $I_1=\{1,3\}$. 
Then $r_0=1$, $r_1=2$, and
$$
\align
\jj^0 &= (2,1,3,2,1,3), \\
\jj^1 &= (1,3,2,1,3,2), \\
f_{\jj^0}(a_1,a_2,a_3,a_4,a_5,a_6) &= 
\left[
\matrix
1 & 0 & 0 & 0 \\
0 & 1 & a_1 & 0 \\
0 & 0 & 1 & 0 \\
0 & 0 & 0 & 1
\endmatrix
\right]
\left[
\matrix
1 & a_2 & 0 & 0 \\
0 & 1 & 0 & 0 \\
0 & 0 & 1 & 0 \\
0 & 0 & 0 & 1
\endmatrix
\right]
\cdots
\left[
\matrix
1 & 0 & 0 & 0 \\
0 & 1 & 0 & 0 \\
0 & 0 & 1 & a_6 \\
0 & 0 & 0 & 1
\endmatrix
\right] \\
&=\left[
\matrix
1 & a_2+a_5 & a_2a_4 & a_2a_4a_6 \\
0 & 1 & a_1+a_4 & a_1a_3+a_1a_6+a_4a_6 \\
0 & 0 & 1 & a_3+a_6 \\
0 & 0 & 0 & 1
\endmatrix
\right]\!,
\tag 0.4.1 \\
f_{\jj^1}(b_1,b_2,b_3,b_4,b_5,b_6) &= \left[
\matrix
1 & b_1 & 0 & 0 \\
0 & 1 & 0 & 0 \\
0 & 0 & 1 & 0 \\
0 & 0 & 0 & 1
\endmatrix
\right]
\left[
\matrix
1 & 0 & 0 & 0 \\
0 & 1 & 0 & 0 \\
0 & 0 & 1 & b_2 \\
0 & 0 & 0 & 1
\endmatrix
\right]
\cdots
\left[
\matrix
1 & 0 & 0 & 0 \\
0 & 1 & b_6 & 0 \\
0 & 0 & 1 & 0 \\
0 & 0 & 0 & 1
\endmatrix
\right] \\
&=\left[
\matrix
1 & b_1+b_4 & b_1b_3+b_1b_6+b_4b_6 & b_1b_3b_5 \\
0 & 1 & b_3+b_6 & b_3b_5 \\
0 & 0 & 1 & b_2+b_5 \\
0 & 0 & 0 & 1
\endmatrix
\right].
\tag 0.4.2
\endalign
$$

Proposition~0.2 asserts that each of the 6 parameters $a_1,a_2,a_3,a_4,a_5,a_6$ 
(resp., $b_1,b_2,b_3,b_4,b_5,b_6$) can be expressed
as a rational function in the 6 matrix entries $u_{ij}$ ($1\le i<j\le 4$) 
of the unipotent upper-triangular matrix 
$$u=(u_{ij})=f_{\jj^0}(a_1,a_2,a_3,a_4,a_5,a_6)$$
(resp., $f_{\jj^1}(b_1,b_2,b_3,b_4,b_5,b_6)$). 
For example,
$$
\aligned
&a_1=\frac{u_{13}u_{24}-u_{14}u_{23}}{u_{13}u_{34}-u_{14}}\,,\,
a_2=\frac{u_{13}u_{34}-u_{14}}{u_{23}u_{34}-u_{24}}\,,\,
a_3=\frac{u_{13}u_{34}-u_{14}}{u_{13}}\,, \\
&a_4=\frac{u_{13}(u_{23}u_{34}-u_{24})}{u_{13}u_{34}-u_{14}}\,,\,
a_5=u_{12}-\frac{u_{13}u_{34}-u_{14}}{u_{23}u_{34}-u_{24}}\,,\,
a_6=\frac{u_{14}}{u_{13}}\,.
\endaligned
\tag 0.4.3
$$
(For explicit formulas for matrices of arbitrary size, see \cite{BFZ96, Theorem~1.4}.) 

Any rational function 
$$\ph(u)=\ph(u_{12}, u_{13}, u_{14}, u_{23}, u_{24}, u_{34})\in [O(U^+)]
$$
can be rewritten in terms of the parameters~$a_i$ (resp.,~$b_i$),
by substituting the appropriate expressions for the $u_{ij}$ from (0.4.1)--(0.4.2):
$$ \aligned
\ph(u)&=\ph(a_2+a_5, \, a_2a_4, \, a_2a_4a_6, \, a_1+a_4, \, a_1a_3+a_1a_6+a_4a_6, \, a_3+a_6)\\
&=\ph(b_1+b_4 ,\, b_1b_3+b_1b_6+b_4b_6 ,\, b_1b_3b_5,\, b_3+b_6 ,\, b_3b_5,\, b_2+b_5). 
\endaligned
\tag 0.4.4$$ 
Theorem~0.3 asserts that $\phi$ is a polynomial in the variables $u_{ij}$
if and only if both functions in the parameters~$a_i$ (resp.,~$b_i$) appearing in (0.4.4) are polynomial. 

Theorem~0.3 can also be restated entirely in terms of the parameters $a_i$ and~$b_i$. 
As observed in \cite{L94} (in a more general setting of an arbitrary
pair of reduced expressions), the birational map relating the
$\nu$-tuples $(a_i)$ and $(b_j)$ to each other 
can be obtained as a composition of simple birational transformations associated to individual braid moves. 
In our example, calculations based on those rules yield the following formulas
expressing $a_1,a_2,a_3,a_4,a_5,a_6$ in terms of $b_1,b_2,b_3,b_4,b_5,b_6$:
$$
a_1= \frac{b_3 b_4 b_5 b_6}{R}, 
a_2= \frac{R}{Q}, 
a_3= \frac{R}{P}  , 
a_4= \frac{PQ}{R} , 
a_5=  \frac{b_2b_3b_4}{Q}, 
a_6=  \frac{b_1b_3b_5}{P}, 
\tag 0.4.5
$$
where
$$
\aligned
P&=b_1 b_3+b_1 b_6+ b_4 b_6 \,,\\
Q&=b_2 b_3+b_2 b_6+ b_5 b_6 \,,\\
R&=b_1 b_2 b_3+b_1b_2b_6+b_1b_5b_6+b_2b_4b_6+b_4b_5b_6
\endaligned
\tag 0.4.6
$$
Theorem~0.3 (in this example) says
that a polynomial $\Phi(a_1,a_2,a_3,a_4,a_5,a_6)$
lies in the subring 
$$
\CC[a_2+a_5, a_2a_4, a_2a_4a_6, a_1+a_4, a_1a_3+a_1a_6+a_4a_6, a_3+a_6]
\subset \CC[a_1,a_2,\dots,a_6]
$$
(cf.\ (0.4.1))
if and only if substituting (0.4.5)--(0.4.6) into $\Phi(a_1,a_2,a_3,a_4,a_5,a_6)$ 
produces a {\it polynomial}
(rather that merely a rational function) 
in $b_1,b_2,\dots,b_6$. 
(An alternative criterion would be to substitute (0.4.3) into $\Phi(a_1,a_2,a_3,a_4,a_5,a_6)$
and verify that the result lies in $\CC[u_1,u_2,\dots,u_6]$.)

\enddefinition


\head 1. Preliminaries on the Weyl group and weights\endhead
\subhead 1.1\endsubhead
Let $\io:G@>>>G$ be the unique automorphism of $G$ such that
$\io(x_i(a))=y_i(a)$,  $\io(y_i(a))=x_i(a)$ for $i\in I,a\in\CC$ and
$\io(t)=t\i$ for $t\in T$. We have $\io^2=1$.

Let $\cy=\Hom(\CC^*,T)$ and $\cx=\Hom(T,\CC^*)$. 
We write the operation in each of these
groups as addition. Let $\la\,,\ra:\cy\T\cx@>>>\ZZ$ be the obvious perfect
pairing.
For $i\in I$, let $\a_i\in\cx$ be the simple root corresponding to $U_i$ and let
$\cha_i$ be the corresponding simple coroot.
Let $\cx^+=\{\l\in\cx\mid \la\cha_i,\l\ra\ge0\ \frl i\in I\}$.
For $i\in I$, the fundamental weight $\o_i\in\cx$ is defined by the condition $\la\cha_j,\o_i\ra=\d_{ij}$ for
$j\in I$. We have $\o_i\in\cx^+$.

For $i\in I$, we denote by $P_i$ the (parabolic) subgroup of $G$ 
generated by $B^+$ together with $\bigcup_{j\in I-\{i\}}U^-_j$.

\subhead 1.2\endsubhead
For $i\in I$ define $s_i:\cy@>>>\cy$ by $\c\m\c-\la\c,\a_i\ra\cha_i$. 
Let $W$ be the subgroup of $\operatorname{Aut}(\cy)$ generated
by $\{s_i;i\in I\}$. This is a Coxeter group with the simple reflections
$\{s_i\mid i\in I\}$ and with the length function that we denote by $w\m|w|$.
Let $w_0\in W$ be the unique element such that $|w_0|=\nu$,
the maximal possible length.
Now $W$ acts on $\cx$ by
the rule $\la\c,w(\l)\ra=\la w\i(\c),\l\ra$ for $\c\in\cy,\l\in\cx$.
For $i\in I$, let $W\o_i$ be the $W$-orbit of the weight $\o_i$ in $\cx$ and let
$W_{I-\{i\}}$ be the subgroup of $W$ generated by
$\{s_j;j\in I-\{i\}\}$.
This is exactly the stabilizer of $\o_i$ with respect to the $W$-action on $\cx$.

Let $NT$ be the normalizer of $T$ in $G$. Now $NT/T$ acts in an obvious way on $\cy$.
This gives an embedding of $NT/T\hookrightarrow\operatorname{Aut}(\cy)$ that identifies $NT/T$ with~$W$. 

For~$i\in I$, set
$\ds_i=x_i(1)y_i(-1)x_i(1)\in NT$ and 
$\ddot s_i=y_i(1)x_i(-1)y_i(1)\in NT$. 
We extend this to define representatives
$\dw\in NT$ and $\ddot w\in NT$ for all $w\in W$ by requiring that
for any $w,w',w''\in W$ satisfying $w''=ww'$ and $|w''|=|w|+|w'|$,
we have $\dw''=\dw\dw'$ and $\ddot w''=\ddot w \ddot w'$.

For $\e\in\{0,1\}$, we set 
$$
z^\e=\displaystyle\prod_{i\in I_{\e}}s_i\in W.
$$ 
(Here the factors commute, so the order does not matter.)

\proclaim{Lemma 1.3 {\rm (see \cite{Bo, Chapter V, \S6, Ex.~2})}}
We have
$$
|z^{[\e]}z^{[\e+1]}\do z^{[\e+h-1]}|=|z^{[\e]}|+|z^{[\e+1]}|+\do+
|z^{[\e+h-1]}|=\nu. 
\tag 1.3.1
$$
It follows that, if $I'\sub I_{[\e]},
I''\sub I_{[\e+l+1]}$, $w'=\displaystyle\prod_{i\in I'}s_i$,
$w''=\displaystyle\prod_{i\in I''}s_i$, then 
$$
|wz^{[\e+1]}\do z^{[\e+l]}w'|=|w|+|z^{[\e+1]}|+|z^{[\e+2]}|+\do+|z^{[\e+l]}|+|w'|, 
\tag 1.3.2
$$
provided either {\rm (a)} $l\in\{0,1,\do,h-2\}$ or
{\rm (b)} $w=1$ and $l\in\{0,1,\do,h-1\}$, or
{\rm (c)} $w'=1$ and $l\in\{0,1,\do,h-1\}$.
\endproclaim

\subhead 1.4\endsubhead
We denote
$$\align
Y'&=\{\o_i\mid i\in I\}, \tag 1.4.1 \\
Y''&=\{w_0\o_i\mid i\in I\}. \tag 1.4.2 
\endalign$$
If $\g\in Y'$,  then $\la\cha_j,\g\ra\ge0$ for all $j\in I$.
If $\g\in Y''$,  then $\la\cha_j,\g\ra\le0$ for all $j\in I$.

\subhead 1.5\endsubhead
Fix $\e\in\{0,1\}$. 
Recall that $\jj^\e=(j^\e_1,j^\e_2,\do,j^\e_\nu)$ was defined in~(0.1.1).
For $k\in\{1,2,\do,\nu\}$, we set
$$
\align
\g^\e_k&=s_{j_\nu^\e}\do s_{j_{k+1}^\e}s_{j_k^\e}\o_{j_k^\e}, \\
\ti\g^\e_k&=s_{j_\nu^\e}\do s_{j_{k+2}^\e}s_{j_{k+1}^\e}\o_{j_k^\e}.
\endalign$$
In order to represent $\g^\e_k$ and $\ti\g^\e_k$ more explicitly, 
we will need to introduce some additional notation.
For $l\in\{1,2,\do,h-2\}$ and $i\in I_{[\e+h-l]}$, let 
$$\align
v_{l,i}^\e&=z^{[\e+h-1]}z^{[\e+h-2]}\cdots z^{[\e+h-l+1]}s_i\o_i.
\endalign$$
Let $\cx_1^\e\sqc\cx_2^\e\do\sqc\cx_h^\e$ be the partition of $\{1,2,\do,\nu\}$ given by
$$\align
\cx_1^\e&=\{1,2,\do,r_{[\e]}\},\\
\cx_2^\e&=\{r_{[\e]}+1,r_{[\e]}+2,\do,r_{[\e]}+r_{[\e+1]}\},\\ 
\cx_3^\e&=\{r_{[\e]}+r_{[\e+1]}+1,r_{[\e]}+r_{[\e+1]}+2,\do,r_{[\e]}+r_{[\e+1]}+r_{[\e+2]}\}, \\
\cdots & \cdots\cdots\cdots
\endalign
$$
Since $s_js_{j'}=s_{j'}s_j$ for $j,j'$ in the same $I_{\e}$ and
$s_j\o_{j'}=\o_{j'}$ if $j\ne j'$, we see that 
$$\align
\g^\e_k=v_{l,j_k^\e}^\e  &\ \ \text{if}\ \  l\in\{1,2,\do,h-2\}, \ k\in\cx_{l+2}^\e\sub\{r+1,r+2,\do,\nu\}, \\
\ti\g^\e_k=v_{l,j_k^\e}^\e &\ \ \text{if}\ \  l\in\{1,2,\do,h-2\}, \ k\in\cx_l^\e\sub\{1,2,,\do,\nu-h\}\\ 
\ti\g^\e_k\in Y'&\ \ \text{if}\ \  k\in\cx_{h-1}^\e\sqc\cx_h^\e=\{\nu-h+1,\nu-h+2,\do,\nu\},\\ 
\g^\e_k\in Y''&\ \ \text{if}\ \  k\in\cx_1^\e\sqc\cx_2^\e=\{1,2,\do,r\}.
\endalign
$$

For $k,k'$ in $\{1,2,\do,\nu\}$ such that $(j_k^\e,j_{k'}^\e)\n I^*$
(see Definition~0.1), we set 
$$\g^\e_{k,k'}=s_{j_\nu^\e}\do s_{j_{k+1}^\e}s_{j_k^\e}\o_{j_{k'}^\e}.$$
From the definitions we see that under these assumptions,

\noindent
(a) $\g^\e_{k,k'}$ is either equal to one of the elements $\g^\e_{k''}$ or lies in $Y'$.

\proclaim{Lemma 1.6} Let $\g=v_{l,i}^\e$ where $\e\in\{0,1\}$,
$i\in I_{[\e+h-l]}$, $l\in\{1,2,\do,h-2\}$.

\noindent
{\rm (a)} If $j\in I_{[\e+h]}$, then $\la\cha_j,\g\ra\ge0$.

\noindent
{\rm (b)} There exists $j\in I_{[\e+h+1]}$ such that $\la\cha_j,\g\ra<0$.
\endproclaim

\demo{Proof}
Let us prove (a). We have
$$\la\cha_j,\g\ra=\la s_iz^{[\e+h-l+1]}\do z^{[\e+h-1]}\cha_j,\o_i\ra.
$$
To show that this is nonnegative, it suffices to prove that

\noindent
{\rm (c)} $s_iz^{[\e+h-l+1]}\do z^{[\e+h-1]}\cha_j$ is a positive coroot.
\nl
We have
$$|s_iz^{[\e+h-l+1]}\cdots z^{[\e+h-1]}|=|s_i|+|z^{[\e+h-l+1]}|+\do+
|z^{[\e+h-1]}|.
$$
(Use $i\in I_{[\e+j-l]}$ and (1.3.2) which holds since $l-1\le h-1$.)
Therefore, to prove~(c), it is enough to show that 
$$|s_iz^{[\e+h-l+1]}\cdots z^{[\e+h-1]}s_j|
=|s_i|+|z^{[\e+h-l+1]}|+\do+|z^{[\e+h-1]}|+|s_j|.
$$
The latter follows from (1.3.2) since $l-1\le h-2$. This proves (a).

Now suppose that (b) does not hold. 
Then by~(a), we have
$\la\cha_j,\g\ra\ge0$ for every $j\in I$. Therefore $\g\in\cx^+$. Since
$\g\in W\o_i$, we have $\g=\o_i$. Hence
$z^{[\e+h-1]}\cdots z^{[\e+h-l+1]}s_i$ is in the stabilizer of $\o_i$, i.e., in $W_{I-\{i\}}$. 
This contradicts
$$|z^{[\e+h-1]}\cdots z^{[\e+h-l+1]}s_i|=|z^{[\e+h-1]}|+\do+|z^{[\e+h-l+1]}|+|s_i|$$
which holds by (1.3.2). 
\qed 
\enddemo

\proclaim{Lemma 1.7} Let $\e\in\{0,1\}$, $i\in I_{[\e+h-l]}$,
$l\in\{1,2,\do,h-2\}$. Let $w\in W$ be the unique element of minimal
length in $\{w_1\in W\mid w_1\o_i=v_{l,i}^\e\}$.

\noindent
{\rm (a)} If $j\in I_{[\e+h]}$, then $|s_jw|>|w|$.

\noindent
{\rm (b)} There exists $j\in I_{[\e+h+1]}$ such that $|s_jw|<|w|$.
\endproclaim

\demo{Proof}
Assume that $j\in I$ satisfies $|s_jw|<|w|$. Then $|w\i s_j|<|w\i|$, 
and using \cite{BZ97, Proposition~2.6} we see that $\la\cha_j,v_{l,i}^\e\ra<0$.
Now using Lemma~1.6(a), we deduce that $j\n I_{[\e+h]}$, proving~(a).
Now suppose (b) does not hold. 
Then by (a),  $|s_jw|>|w|$ for all
$j\in I$. Hence $w=1$ and $v_{l,i}^\e=\o_i$ so that
$\la\cha_j,v_{\l,i}^\e\ra\ge0$ for all $j\in I$. This contradicts
Lemma 1.6(b). 
\qed\enddemo

\subhead 1.8\endsubhead
Let $\e\in\{0,1\}$. Denote
$$Y^\e=\{v_{l,i}^\e\mid i\in I_{[\e+h-l]},l\in\{1,2,\do,h-2\}\}. \tag 1.8.1
$$
We are going to show that 

\noindent
(a) all the weights in $Y^\e$ are distinct.
\nl
To prove this, suppose that $v_{l,i}^\e=v_{l',i'}^\e$ where 
$i\in I_{[\e+h-l]}$, $i'\in I_{[\e+h-l']}$, and $l,l'\in\{1,2,\do,h-2\}$.
Then $W\o_i=W\o_{i'}$ and therefore $\o_i=\o_{i'}$ and so $i=i'$. 

Suppose that $l\neq l'$. Without loss of generality, we may assume that $l>l'$. 
Setting $e=l-l'\ge1$ we get: 
$$z^{[\e+h-l+e]}z^{[\e+h-l+e-1]}\cdots z^{[\e+h-l+1]}s_i\o_i=s_i\o_i .
$$
Hence

\noindent
(b) $s_iz^{[\e+h-l+e]}z^{[\e+h-l+e-1]}\do z^{[\e+h-l+1]}s_i\in
W_{I-\{i\}}$.
\nl
From (1.3.2) we see that
$$\align
&|s_iz^{[\e+h-l+e]}z^{[\e+h-l+e-1]}\do z^{[\e+h-l+1]}s_i| \\
=&|s_iz^{[\e+h-l+e]}|+|z^{[\e+h-l+e-1]}|+\do+|z^{[\e+h-l+1]}|+|s_i|
\endalign
$$
which contradicts (b). 
Hence $l=l'$ and (a) is proved. 

We note that

$$\sha(Y^\e)=
\sum_{l=1}^{h-2}r_{[\e+h-l]}=\sum_{l=1}^hr_{[\e+h-l]}-r_{[\e+1]}-r_{[\e]}=\nu-r.
\tag 1.8.2
$$

\proclaim{Lemma 1.9}
With the notation introduced in {\rm (1.4.1)}, {\rm (1.4.2)}, {\rm (1.8.1)}, we have 
$Y^\e\cap Y'=\emp$ and $Y^\e\cap Y''=\emp$ (assuming $G$ is not of type~$A_1$). 
\endproclaim

\demo{Proof}
If $\g\in Y^\e$, then $\la\che\a_j,\gamma\ra<0$ for some $j$ by Lemma~1.6(b);
thus $\g\notin Y'$ by~1.4. 

Assume that $v_{l,i}^\e=w_0\o_j$  for some
$l\in\{1,2,\do,h-2\}$, $i\in I_{[\e+h-l]}$, $j\in I$.
Then $\o_i$ and $\o_j$ are in the same $W$-orbit.
Hence $i=j$ and we have
$$\align
& z^{[\e+h-1]}z^{[\e+h-2]}\cdots z^{[\e+h-l+1]}s_i\o_i \\
=
& z^{[\e+h-1]}z^{[\e+h-2]}\cdots z^{[\e+h-l+1]}z^{[\e+h-l]}\cdots  z^{[\e]}\o_i \,. 
\endalign
$$
This implies that  $s_iz^{[\e+h-l]}\cdots  z^{[\e]}\o_i=\o_i$, i.e., 
$s_iz^{[\e+h-l]}\cdots  z^{[\e]}$ lies is in the stabilizer of~$\o_i$, that is, in $W_{I-\{i\}}$.  
So any reduced expression of it does not contain~$s_i$.
If $l\ge2$, this contradicts
$$  |s_iz^{[\e+h-l]}\cdots  z^{[\e]}|=|s_i|+|z^{[\e+h-l]}|+\cdots +|z^{[\e]}|.
$$
since $h-l+1+1\le h$.  Therefore $l=1$ and moreover
any reduced expression of $s_iw_0$ does not contain~$s_i$.
But this cannot happen if $G$ is of type other than~$A_1$. 
Indeed, for some $\e\in\{0,1\}$, 
$$ z^{[\e+h-1]}z^{[\e+h-2]}\cdots z^{[\e+h-l+1]}z^{[\e+h-l]}\cdots  z^{[\e]}
$$
gives a reduced expression of $w_0$ such that $s_i$ appears in the first
group $z^{[\e+h-1]}$. 
If $s_i$ does not appear in any other group, then 
there are only two factors and $h=2$. But $h>2$ in any type other than~$A_1$.
\qed\enddemo

\head 2. An irreducibility property\endhead

In this section, we prove the following result.

\proclaim{Proposition 2.1} Let $w, w'\!\in\! W$. 
The set $U^+\cap(B^-\dw'B^+\dw\i)$ is empty if $w'\not\le w$; 
it is smooth and irreducible, of dimension $\nu-|w'|$, if $w'\le w$.     
\endproclaim

\subhead 2.2\endsubhead 
For $y\in W$, let
$$\align
U^+_y&=\{u\in U^+,\dy\i u\dy\in U^-\}, \\
U^{+y}&=\{u\in U^+, \dy\i u\dy\in U^+\}.
\endalign
$$
The multiplication map 
$U^+_y\T U^{+y}@>\si>>U^+$
is an isomorphism of varieties.

\subhead 2.3\endsubhead
For $x\in G$ and a subgroup $C$ of $G$, we shall write ${}^xC$ instead of
$xCx\i$. For $w\in W$, we shall write ${}^wC$ instead of ${}^{\dw}C$.

We denote by $\cb$ the variety of Borel subgroups in $G$. For $B',B''\in\cb$, 
there is a unique $w\in W$ such that for some $x',x''$ in $G$ we have
$B'={}^{x'}B^+, B''={}^{x''}B^+$, $x'{}\i x''\in B^+\dw B^+$; we then write
$w=\operatorname{pos}(B',B'')$. 

For $z,z'\in W$, we denote
$$
\car_{z,z'}=\{B\in\cb \mid \operatorname{pos}(B^-,B)=z',
\ \operatorname{pos}(B,B^+)=z\i w_0\}.
$$
It is known \cite{KL79} that $\car_{z,z'}$ is nonempty if and only if
$z\le z'$. We show:

\proclaim{Proposition 2.4}
If $z\le z'$ then $\car_{z,z'}$ is smooth, irreducible of dimension
$|z'|-|z|$.
\endproclaim

\demo{Proof}
We shall adapt an argument in \cite{L98, 1.4} by replacing $\RR$ by $\CC$. The set 
$$
\ti\car_{z,z'}=\{B\in\cb\mid \operatorname{pos}(B^-,B)=z',\operatorname{pos}(B,{}^{w_0z}B^-)=w_0\}.
$$
is an open nonempty subset in $\{B\mid \operatorname{pos}(B^-,B)\!=\!z'\}\cong\CC^{|z'|}$.
Hence it is smooth irreducible of dimension $|z'|$. Clearly, the
map $(B,u)\m{}^uB$ is an isomorphism
$\car_{z,z'}\T(U^-\cap{}^{w_0z}U^-)@>\si>>\ti\car_{z,z'}$. 
Now the claim follows since $U^-\cap{}^{w_0z}U^-\cong\CC^{|z|}$.
\qed\enddemo

\subhead 2.5 \endsubhead
A result related to Proposition 2.4 holds for the
analogue of $\car_{z,z'}$ over a finite field~$\FF_q$. By \cite{KL79}, 
the number of $\FF_q$-rational points in this analogue is given by 
the polynomial $R_{z,z'}$ in {\it loc.cit.} evaluated at $q$.
By the inductive
formula in {\it loc.cit.}, the latter is monic of degree $|z'|-|z|$.

\demo{Proof of Proposition~2.1}
Setting $B={}^xB^+$, we can reformulate Proposition 2.4 as the statement that
$$\align
&\{xB^+\in G/B^+ \mid \operatorname{pos}(B-,{}^xB^+)=z',\operatorname{pos}({}^xB^+,B^+)=z\i w_0\}\\&=
((U^+w_0z\dot{}B^+)\cap (B^-(w_0z'{}\i)\dot{}B^+))/B^+\\&
=(U^+_{w_0z}(w_0z)\dot{})\cap(B^-(w_0z'{}\i)\dot{}B^+)\endalign$$
is smooth, irreducible of dimension $|z'|-|z|$ if $z\le z'$, and is
empty if $z\not\le z'$.

Replacing here $w_0z,w_0z'$ by $w,w'$ we deduce that 
$(U^+_w\dw)\cap(B^-\dw'B^+)$ is smooth, irreducible of dimension
$|w|-|w'|$ if $w'\le w$, and is empty if $w'\not\le w$. 

Using 2.2, we see that the map
$$(U^+_w \dw)\cap(B^-\dw'B^+)\T U^{+w}@>>>(U^+\dw)\cap(B^-\dw'B^+)$$
given by $(u'\dw,u'')\m u'u''\dw$ with $u'\in U^+_w$ such that
$u'\dw\in B^-\dw'B^+$ and $u''\in U^{+w}$ is an isomorphism of varieties.
Since $U^{+w}\cong\CC^{\nu-|w|}$, we conclude that
$(U^+\dw)\cap(B^-\dw'B^+)$ is smooth, irreducible of dimension
$\nu-|w'|$ if $w'\le w$, and is empty if $w'\not\le w$. 
This completes the proof of Proposition 2.1.
\qed\enddemo

\head 3. Proof of Theorem 0.3\endhead

When $G$ is of type $A_1$, we have $\jj^0=\jj^1$ and the theorem is trivial.
For the rest of this section, we assume that $G$ is not of type~$A_1$. 

\subhead 3.1\endsubhead
Fix $i\in I$. Let 
$$V_i=\{f\in O(G) \mid f(utg)=\o_i(t)f(g) \ \frl u\in U^-,t\in T,g\in G\}.$$
The group $G$ acts on $V_i$ by $g_1:f\m g_1f$ where $(g_1f)(g)=f(gg_1)$.
There is a unique $f\in V_i$ such that $f(gu)=f(g)$ for all
$g\in G,u\in U^+$ and such that $f(1)=1$. We denote it by~$\D$.
(Note that $\D$ depends on the choice of~$i$.)

We show that $\D(\ds_i)=0$. Setting $g_c=y_i(-c)\cha_i(c\i)x_i(c)$
for $c\in\CC^*$, we~see that $\displaystyle\lim_{c\to\iy}g_c=\ds_i$ in $G$. We have
$\D(g_c)\!=\!\o_i(\cha_i(c\i))\!=\!c\i$, so $\D(\ds_i)\!=\!\displaystyle\lim_{c\to\iy}c\i\!=\!0$.
It follows that $\D$ vanishes on $U^-\ds_iB^+$, hence also on the closure

\noindent
(a) $Z=\ov{U^-\ds_iB^+}=\cup_{w;s_i\le w}U^-\dw B^+
=\cup_{w\in W-W_{I-\{i\}}}U^-dw B^+=G-(U^-P_i)$.
\nl
The function $\D$ is preserved (up to a nonzero scalar) by the action
of $P_i$ on $V_i$. Hence $\D$ takes only nonzero values on the open subset
$U^-P_i$ of $G$, implying that

\noindent
(b) $Z=\{g\in G;\D(g)=0\}$.

\definition{Definition 3.2}
Let $i\in I$ and $\g\in W\o_i$. 
Following \cite{BZ97}, we set $\D_\g=\ddot w\D$, where $w\in W$
is such that $w\o_i=\g$. This does not depend on the choice of $w$.
In particular, $\D_{\o_i}=\D$.  

Let $\D^+_\g$ be the restriction of $\D_\g$ to $U^+$. For $u\in U^+$, we
have $\D^+_\g(u)=\D(u\ddot w)$, with $w$ as above.
(Note that $\D^+_\g$ is not identically zero on $U^+$.
Otherwise we would have $\D(U^-B^+\ddot w)=0$; but $U^-B^+\ddot w$
is dense in $G$; hence $\D=0$, a contradiction.) 

We will also use the notation
$$\cz_\g=\{u\in U^+\mid \D^+_\g(u)=0\}. \tag 3.2.1
$$
\vskip -.1in
\enddefinition

\proclaim{Lemma 3.3}
Let $i\in I$, $w\in W$, and $\g=w\o_i$. Then

\noindent
{\rm (a)} $\cz_\g=\bigcup_{y\in W-W_{I-\{i\}}}(U^+\cap(U^-\dy B^+\dw\i))$;


\noindent
{\rm (b)} if $s_i\not\le w$ then $\cz_\g$ is empty; 

\noindent
{\rm (c)} if $s_i\le w$, then $\cz_\g$ is the closure of
$U^+\cap(U^-\ds_iB^+\dw\i)$ (a smooth irreducible variety of dimension $\nu-1$). 
\endproclaim

\demo{Proof}
Using 3.1(a),(b), we get
$$\align
\cz_\g&=\{u\in U^+;\D(u\ddot w)=0\} \\
&=\{u\in U^+;u\ddot w\in Z\}\\
&= \{u\in U^+\mid u\ddot w\in\cup_{y\in W-W_{I-\{i\}}}U^-dy B^+\}, 
\endalign$$
and (a) follows. (We used that $B^+\ddot w\i=B^+\dot w\i$.)

By Proposition 2.1, $U^+\cap (U^-\ds_i B^+\dw\i)$ is smooth irreducible
of dimension $\nu-1$ provided that $s_i\le w$ and is empty if
$s_i\not\le w$. Moreover if $y$ satisfies $s_i<y$, then the same
Proposition shows that $U^+\cap(U^-\dy B^+\dw\i)$ is either empty or
irreducible of dimension $\nu-|y|\le\nu-2$.
Since, by Krull's theorem, $\cz_\g$ is either empty or of pure dimension
$\nu-1$, the statements (b) and (c) follow. 
\qed\enddemo

\proclaim{Lemma 3.4} Let $\e\in\{0,1\}$, $\g\in Y^\e$. Then: 

\noindent
{\rm (a)} $\cz_\g$ (see {\rm(3.2.1)}) is an irreducible variety of dimension $\nu-1$;

\noindent
{\rm (b)} for any $j\in I_{[\e+h]}$ and any $c\in\CC$ we have
$\cz_\g x_j(c)\sub\cz_\g$;

\noindent
{\rm (c)} there exists $j\in I_{[\e+h+1]}$ such that for some $c\in\CC$ we have
$\cz_\g x_j(c)\not\sub\cz_\g$.
\endproclaim

\demo{Proof of {\rm(a)}}
We write $\g=w\o_i$ with $i\in I,w\in W$. By Lemma~1.9, we have
$\g\n Y'$
hence $w\n W_{I-\{i\}}$ so that $s_i\le w$. Now (a) follows from Lemma 3.3(c).
\qed
\enddemo

\demo{Proof of {\rm(b)}}
We write $\g=w\o_i$ where $i\in I$ and $w\in W$ is the
unique element of minimal length in $\{w_1\in W\mid w_1\o_i=\g\}$. Using
Lemma 3.3(c), we see that it is enough to show that for $j,c$ as in (b) we have
$$(U^+\cap(U^-\dy B^+\dw\i))x_j(c)\sub U^+\cap(U^-\dy B^+\dw\i)$$
for any $y\in W-W_{I-\{i\}}$.
This follows from $\dw\i x_j(c)\in U^+\dw\i$ which in turn follows from
$|s_jw|>|w|$ (see 1.7(a)) or equivalently $|w\i s_j|>|w\i|$. 
\qed
\enddemo

\demo{Proof of {\rm(c)}}
Suppose that (c) does not hold. Using (b), we see that
for any $j\in I$ and any $c\in\CC$ we have $\cz_\g x_j(c)\sub\cz_\g$.
Since the elements $x_j(c)$ for various $j,c$ generate the group $U^+$, it
follows that $\cz_\g U^+\sub\cz_\g$. Since $\cz_\g\ne\emp$, we conclude that
$\cz_\g=U^+$.
This contradicts Lemma 3.3(b),(c). 
\qed
\enddemo

\proclaim{Lemma 3.5}
Let $\g\in Y^0$ and $\g'\in Y^1$. Then every irreducible
component of $\cz_\g\cap\cz_{\g'}$ has dimension $\le\nu-2$.
\endproclaim

\demo{Proof}
By Lemma 3.4(c) with $\e=0$, there exist $j\in I_{[h+1]}$ and $c\in\CC$ such that
$\cz_\g x_j(c)\not\sub\cz_\g$. By  Lemma 3.4(b) with $\e=1$, we have 
$\cz_{\g'} x_j(c)\sub\cz_{\g'}$. Therefore $\cz_\g\ne\cz_{\g'}$.
Since $\cz_\g,\cz_{\g'}$ are irreducible of dimension $\nu-1$, the lemma
follows.
\qed\enddemo

\subhead 3.6 \endsubhead
Consider the partition
$$
U^+=\bigsqcup_{z\in W}U^+(z)
\tag 3.6.1
$$
where $$U^+(z)=U^+\cap B^-\dz B^-$$ is smooth and irreducible of dimension
$|z|$ (cf.\ Proposition 2.1 with $(w,w')$ replaced by $(w_0,zw_0)$).
Furthermore, the closure of $U^+(z)$ in $U^+$
is equal to $\displaystyle\bigsqcup_{z';z'\le z}U^+(z')$. It follows that $U^+(w_0)$ is open
dense in $U^+$.
For $z\in W$, we set $$U^-(z)=U^-\cap B^+\dz B^+=\io(U^+(z))$$
(see~1.1 for the definition of~$\io$).
Then $$U^-=\displaystyle\bigsqcup_{z\in W}U^-(z)$$ and $U^-(w_0)$ is open dense in $U^-$.

Let
$$A:U^+(w_0)@>\si>>U^+(w_0)$$
be the composition
$$U^+(w_0)@>\si>>
\{B\in\cb \mid \operatorname{pos}(B^+,B)=
\operatorname{pos}(B,B^-)=w_0\}@>\si>>U^-(w_0)@>\si>>U^+(w_0)$$
where
the first isomorphism is $u\m{}^uB^-$, 
the second isomorphism is the inverse of $u'\m{}^{u'}B^+$,
and the third isomorphism is the restriction of $\io$. 

We will show that $A$ is an involution.
For $u\in U^+$, we have ${}^uB^-={}^{\io(A(u))}B^+$. 
Replacing $u$ by $A(u)$ we
obtain ${}^{A(u)}B^-={}^{\io(A^2(u))}B^+$. Applying $\io$, we obtain
${}^{\io(A(u))}B^+={}^{A^2(u)}B^-$, i.e., ${}^uB^-={}^{A^2(u)}B^-$.
Hence $u=A^2(u)$ and $A^2=1$. 

\subhead 3.7\endsubhead
Let $\e\in\{0,1\}$. We denote
$$\cv^\e=\{u\in U^+ \mid
\D^+_{\g^\e_k}(u)\ne0,\D^+_{\ti\g^\e_k}(u)\ne0,k=\{1,2,\do,\nu\}\}.$$
This set is open in $U^+$. It is also nonempty, since each of
$\D^+_{\g^\e_k},\D^+_{\ti\g^\e_k}$ is not identically zero on $U^+$.
We denote
$$\cv^\e_*=\{u\in U^+(w_0) \mid A(u)\in\cv^\e\}=U^+(w_0)\cap A\i(\cv^\e)\}.$$
This set is open in $U^+$. It is also  nonempty, as it is the
intersection of two open nonempty subsets of $U^+$.
We shall need the following result from \cite{BZ97}, \cite{FZ99}:

\proclaim{Lemma 3.8}
The map $f_{\jj^\e}:\CC^\nu@>>>U^+$ restricts to an isomorphism
$(\CC^*)^\nu@>\si>>\cv^\e_*$.
\endproclaim

\subhead 3.9 \endsubhead
Using the results in 1.5, we see that
$$\cv^\e_*=\{u\in U^+(w_0) \mid \D^+_\g(Au)\ne0 \text{ for all }
\g\in Y^\e\cup Y'\cup Y''\}.\tag 3.9.1$$
If $\g=\o_i$ with $i\in I$, then $\D^+_\g$ is the function
$u\m\D_{\o_i}(u)=\D_{\o_i}(1)=1$ (a~constant function).
If $\g=w_0\o_i$ with $i\in I$, then $\D^+_\g(u)\ne0$ for any
$u\in U^+(w_0)$. (Indeed, writing $u=u'\dw_0b'$ with
$u'\in U^-,b'\in B^-$, so that $u\dw_0=u'tu_1$ with $t\in T,u_1\in U^+$,
we have $\D^+_\g(u)=\D_{\o_i}(u\dw_0)=\D_{\o_i}(u'tu_1)=\o_i(t)\ne0$.)
It follows that $Y'$ and $Y''$ can be eliminated from (3.9.1), and we conclude that 
$$\cv^\e_*=\{u\in U^+(w_0)  \mid \D^+_\g(Au)\ne0 \text{ for all }\g\in Y^\e\}.
\tag 3.9.2$$

\proclaim{Lemma 3.10}
$\dim(U^+(w_0)-(\cv^0_*\cup\cv^1_*))\le\nu-2$.
\endproclaim

\demo{Proof}
From (3.9.2), we obtain
$$U^+(w_0)-(\cv^0_*\cup\cv^1_*)=\bigcup_{(\g,\g')\in Y^0\T Y^1}
A(U^+(w_0)\cap\cz_\g\cap\cz_{\g'}).$$
It remains to use that $\dim(\cz_\g\cap\cz_{\g'})\le\nu - 2$
for $(\g,\g')\in Y^0\T Y^1$ (see~3.5).
\qed\enddemo

\proclaim{Lemma 3.11}
Let $\ii=(i_1,i_2,\do,i_n)\in I^n$ be a reduced expression, that is, the element 
$w=s_{i_1}\do s_{i_n}\in W$
has length $n$. Let 
$${}'f_\ii:(\CC^*)^n@>>>U^+$$ 
be the restriction
of the map $f_\ii$ in {\rm (0.1.2)}. 
Then ${}'f_\ii$ is an isomorphism of $(\CC^*)^n$ onto an open
subset ${}'U^+_\ii$ of $U^+(w)$.
\endproclaim

\demo{Proof}
Induction on $n$. For $n=0$, the result is obvious.
Assume that $n\ge1$. Let $\ii'=(i_1,i_2,\do,i_{n-1})\in I^{n-1}$ and let
$w'=s_{i_1}\do s_{i_{n-1}}$.
The map 
$$\align
U^+(w')\T\CC^* &@>>>U^+(w) \\
(u',c) &\m u'x_{i_n}(c)
\endalign
$$ 
is an isomorphism of $U^+(w')\T\CC^*$ onto an open subset of $U^+(w)$. It
restricts to an isomorphism of ${}'U^+_{\ii'}\T\CC^*$ onto an open subset
${}'U^+_\ii$ of $U^+(w)$. 
\qed\enddemo

\subhead 3.12 \endsubhead
Let $\e\in\{0,1\}$ and let $i\in I_{[\e+h+1]}$. 
Define $k\in\cx_h^{\epsilon}$ (in the notation of 1.5) by $j^\e_k=i$.
Let $\CC^\nu_i$ (resp. ${}'\CC^\nu_i$) be the subset of $\CC^\nu$
consisting of all $(a_1,a_2,\do,a_\nu)$ such that
$a_l\in\CC^*$ for $l\ne k$ whereas $a_k\in\CC$ (resp. $a_k=0$).
By restricting $f_{\jj^\e}:\CC^\nu@>>>U^+$ to
$\CC^\nu_i$ (resp. ${}'\CC^\nu_i$), we obtain maps
$f_{\jj^\e;i}:\CC^\nu_i@>>>U^+$ and ${}'f_{\jj^\e;i}:\CC^\nu_i@>>>U^+$.

It follows from Lemma 3.11 that

\noindent
(a) ${}'f_{\jj^\e;i}$ is an isomorphism of ${}'\CC^\nu_i$ onto an open
subset ${}'U^+_{\jj^\e;i}$ of $U^+(w_0s_i)$.

We next prove that 

\noindent
(b) $f_{\jj^\e;i}$ is an isomorphism of $\CC^\nu_i$ onto an open subset
$U^+_{\jj^\e;i}$ of $U^+(w_0)\cup U^+(w_0s_i)$ containing
${}'U^+_{\jj^\e;i}$.

\demo{Proof}
The map $U^+(w_0s_i)\T\CC@>>>U^+$, $(u',c)\m u'x_{i_n}(c)$, is an
isomorphism of $U^+(w_0s_i)\T\CC$ onto an open subset of
$U^+(w_0)\cup U^+(w_0s_i)$. It
restricts to an isomorphism of ${}'U^+_{\jj^\e;\ii}\T\CC$ onto an open
subset $U^+_{\jj^\e;i}$ of $U^+(w_0)\cup U^+(w_0s_i)$.
\qed\enddemo

The following is a special case of Lemma 3.11: 

\noindent
(c)  ${}'f_{\jj^\e}$ is an isomorphism of $(\CC^*)^\nu$ onto an open
subset ${}'U^+_{\jj^\e}$ of $U^+(w_0)$.

From (c), we deduce that

\noindent
(d) $f_{\jj^\e}$ is a birational isomorphism from $\CC^\nu$ to $U^+$.

\proclaim{Lemma 3.13}
Let \ $\cu$ be the open subset of $U^+$ defined by
$$
\cu=\cv^0_*\cup\cv^1_*\cup\bigcup_{\e\in\{0,1\},i\in I_{[\e+h+1]}}
U^+_{\jj^\e;i}.
\tag 3.13.1
$$
Then \ $\dim(U^+-\cu)\le\nu-2$.
\endproclaim

\demo{Proof}
Using the partition (3.6.1), it is enough to show that 
$$
\dim(U^+(z)\cap(U^+-\cu))\le\nu-2
\tag 3.13.2
$$
for any $z\in W$. 

Case 1: $z=w_0$. We have
$$U^+(w_0)\cap(U^+-\cu)\sub U^+(w_0)-(\cv^0_*\cup\cv^1_*).
$$
Therefore
$$\dim(U^+(w_0)\cap(U^+-\cu))\le\dim(U^+(w_0)-(\cv^0_*\cup\cv^1_*))\le
\nu-2$$
(see Lemma~3.10), and (3.13.2) follows.

Case 2:  $z=w_0s_i$ with $i\in I$.
Define $\e\in\{0,1\}$ by $i\in I_{[\e+h+1]}$. Then
$$U^+(w_0s_i)\cap(U^+-\cu)\sub U^+(w_0s_i)
\cap(U^+-U^+_{\jj^\e;i})\sub U^+(w_0s_i)-{}'U^+_{\jj^\e;i}\,.$$
The last difference has dimension $\le\nu-2$ (as desired) since
$U^+(w_0s_i)$ is irreducible of dimension $\nu-1$ and
${}'U^+_{\jj^\e;i}$ is a nonempty open subset of
$U^+(w_0s_i)$. 

Case 3:   $z$ is not of the form $w_0$ or $w_0s_i$.
Then $|z|\le\nu-2$.
Therefore $\dim(U^+(z))\le \nu-2$ which implies 
(3.13.2).
\qed\enddemo

\demo{{\bf 3.14.} Proof of Theorem 0.3}
The ``only if'' part of Theorem~0.3 is obvious.
Let us prove the ``if'' statement.
Consider $\ph\in[O(U^+)]$ such that
$f_{\jj^\e}^*(\ph)\in[O(\CC^\nu)]$ belongs to $O(\CC^\nu)$ for $\e=0$
and for $\e=1$.
From our assumption we see that $\ph|_{\cv^\e_*}$ is regular for
$\e\in\{0,1\}$ (see Lemma~3.8) and that $\ph|_{U^+_{\jj^\e;i}}$ is regular
for $\e\in\{0,1\}$ and $i\in I_{[\e+h+1]}$ (see 3.12(b)). Hence $\ph$
is regular on $\cu$. Using this and Lemma~3.13, we conclude that $\ph$ is regular
on $U^+$. Theorem~0.3 is proved.
\qed\enddemo

\head 4. The study of $O(G/U^-)$ \endhead
\subhead 4.1\endsubhead
For $\ii=(i_1,i_2,\do,i_\nu)\in I^\nu$, we define the maps 
$$\align
&f_{\ii;+}:\CC^\nu\T T@>>>G/U^- \\
&f_{\ii;-}:\CC^\nu\T T@>>>G/U^-
\endalign
$$
by
$$\align
f_{\ii;+}(a_1,a_2,\do,a_\nu,t)&=
x_{i_1}(a_1)x_{i_2}(a_2)\do x_{i_\nu}(a_\nu) \, t\, U^-,\\
f_{\ii;-}(a_1,a_2,\do,a_\nu,t)&=
y_{i_1}(a_1)y_{i_2}(a_2)\do y_{i_\nu}(a_\nu)\, t\, \dw_0\, U^-.
\endalign
$$
Of particular interest to us are the cases where $\ii=\jj^\e$, for $\e\in\{0,1\}$, 
as in (0.1.1). 
Proposition~0.2 implies that
both $f_{\jj^\e;+}$ and  $f_{\jj^\e;-}$ are birational isomorphisms from
$\CC^\nu\T T$ to $G/U^-$. 
Consequently the maps $f_{\jj^\e;+}^*$ and  $f_{\jj^\e;-}^*$
are well-defined isomorphisms $[O(G/U^-)]@>\si>>[O(\CC^\nu\T T)]$.

\proclaim{Theorem 4.2} 
An element $\ph\in[O(G/U^-)]$ belongs to $O(G/U^-)$ if
and only if each of the four rational functions 
\hbox{$f_{\jj^0;+}^*(\ph), f_{\jj^1;+}^*(\ph), f_{\jj^0;-}^*(\ph), f_{\jj^1;-}^*(\ph)\!\in\![O(\CC^\nu\!\T\! T)]$}
belongs to $O(\CC^\nu\T T)$. 
\endproclaim

The proof of Theorem 4.2 will rely on the following statement. 

\proclaim{Lemma 4.3} We have
$$\dim(G/U^--((U^+TU^-\cup U^-T\dw_0U^-)/U^-))\le
\dim(G/U^-)-2.
\tag 4.3.1$$
\endproclaim

\demo{Proof}
The inequality (4.3.1) is equivalent to 
$$\dim(G/B^--((U^+B^-\cup(U^-\dw_0B^-)/B^-))\le\dim(G/B^-)-2, $$ 
which is equivalent to the inequality 
$$\dim(\cb-(\{B\in\cb\mid \operatorname{pos}(B,B^+)=w_0\}\cup\{B\in\cb;\operatorname{pos}(B,B^-)=w_0\})
\le\dim\cb-2$$
and thus to the statement that, for any $z\in W-\{1\}$ and $z'\in W-\{w_0\}$, we have 
$$\dim(\{B\in\cb\mid \operatorname{pos}(B^-,B)=z',
\ \operatorname{pos}(B,B^+)=z\i w_0\})\le\nu-2.
$$
The last claim follows from Proposition~2.4
since $|z'|-|z|\le\nu-2$. 
\qed
\enddemo

\demo{{\bf 4.4.} Proof of Theorem 4.2}
The ``only if'' statement in the theorem is obvious.
Let us prove the ``if'' statement. 
Thus, let $\ph\in[O(G/U^-)]$ be 
such that the four conditions in the theorem are satisfied.
We need to show that $\ph\in O(G/U^-)$.

Suppose $G$ is of type $A_1$.
Then $\ph$ is regular on
$(U^+TU^-\cup U^-T\dw_0U^-) /U^-$.
Hence by~(4.3.1), it is regular on $G/U^-$, and we are done.

In the rest of the proof, we assume that $G$ is of type other than~$A_1$.

We first show that $\ph$ regular on the open subset
$U^+TU^-/U^-$ of~$G/U^-$. With the notation as in 3.7 and 3.12, we see as in the
proof in 3.14 that $\ph$ is regular on each of the following open
subsets of $U^+TU^-/U^-$: 

$\bullet$\ $\cv^\e TU^-/U^-$, for $\e\in\{0,1\}$; 

$\bullet$\ $U^+_{\jj^\e;i}TU^-/U^-$,
for $\e\in\{0,1\}$ and $i\in I_{[\e+h+1]}$.

\noindent
Hence $\ph$ is regular on the union of these
subsets, i.e., on $\cu T U^-/U^-$ (here $\cu\sub U^+$ is as in (3.13.1)).
By Lemma 3.13, we have 
$$\dim((U^+TU^- - \cu TU^-)/U^-)\le\nu+r-2=\dim(G/U^-)-2.
$$
Since $\ph$ is regular on $\cu TU^-/U^-$, it follows that $\ph$ is
regular on $U^+TU^-/U^-$.

We next show that $\ph$ is regular on the open subset
$U^-T\dw_0U^-/U^-$ of $G/U^-$.

We denote $\cv^{\e-}_*=\io(\cv^\e_*)\subset U^-$ (cf.~3.7)
and $U^-_{\jj^\e;i}=\io(U^+_{\jj^\e;i})\subset U^-$ (cf.~3.12). 


As in the
proof in 3.14 (with $U^+$ replaced by $U^-$), we see that $\ph$ is regular on
each of the following open subsets of $U^-T\dw_0U^-/U^-$:

$\bullet$\ $\cv^{\e-}T\dw_0U^-/U^-$, for $\e\in\{0,1\}$;

$\bullet$\ $U^-_{\jj^\e;i}T\dw_0U^-/U^-$, 
for $\e\in\{0,1\}$ and $i\in I_{[\e+h+1]}$. 

\noindent
Hence $\ph$ is regular on
the union of these subsets, i.e., on $\cu^-T\dw_0U^-/U^-$
where $\cu^-=\io(\cu)\sub U^-$.  By Lemma 3.13 (with $U^-$ instead of $U^+$), 
we have
$$\dim((U^-T\dw_0U^- - \cu^-T\dw_0U^-)/U^-)\le\nu+r-2=\dim(G/U^-)-2.$$
Since $\ph$ is regular on $\cu^-T\dw_0U^-/U^-$, it follows that
$\ph$ is regular on $U^-T\dw_0U^-/U^-$.
Thus $\ph$ is regular on the open subset
$((U^+TU^-)\cup(U^-T\dw_0U^+))/U^-$ of $G/U^-$.
Using this and Lemma 4.3, we conclude that $\ph$ is regular on $G/U^-$, as desired.
\qed\enddemo

\head 5. The study of $O(G)$\endhead
\subhead 5.1\endsubhead
For 
$\ii =(i_1,i_2,\do,i_\nu)\in I^\nu$ and $\ii'=(i'_1,i'_2,\do,i'_\nu)\in I^\nu$,
we define the maps
$$ \align
&f_{\ii,\ii';\pm}:\CC^\nu\T T\T\CC^\nu@>>>G, \\
&f_{\ii,\ii';\mp}:\CC^\nu\T T\T\CC^\nu@>>>G
\endalign
$$
by
$$\align&f_{\ii,\ii';\pm}(a_1,a_2,\do,a_\nu,t,b_1,b_2,\do,b_\nu)\\&=
x_{i_1}(a_1)x_{i_2}(a_2)\do x_{i_\nu}(a_\nu)\, t\,
y_{i'_1}(b_1)y_{i'_2}(b_2)\do y_{i'_\nu}(b_\nu),\\
&f_{\ii,\ii';\mp}(a_1,a_2,\do,a_\nu,t,b_1,b_2,\do,b_\nu)\\&=
y_{i_1}(a_1)y_{i_2}(a_2)\do y_{i_\nu}(a_\nu)\,t\i \,x_{i'_1}(b_1)
x_{i'_2}(b_2)\do x_{i'_\nu}(b_\nu).\endalign$$
Thus, $f_{\ii,\ii';\mp}=\io f_{\ii,\ii';\pm}$.

Let $\jj^\e$, for $\e\in\{0,1\}$, be as in (0.1.1). 
From Proposition~0.2 one can deduce that
for each of the four possible pairs $(\e,\e')\in\{0,1\}\T\{0,1\}$, both maps
$f_{\jj^\e,\jj^{\e'};\pm}$ and  $f_{\jj^\e,\jj^{\e'};\mp}$
are birational isomorphisms from $\CC^\nu\T T\T\CC^\nu$ to $G$.
It follows that both  $f_{\jj^\e,\jj^{\e'};\pm}^*$ and  $f_{\jj^\e,\jj^{\e'};\mp}^*$
are well defined isomorphisms $[O(G)]@>\si>>[O(\CC^\nu\T T\T\CC^\nu)]$.

\proclaim{Theorem 5.2} An element $\ph\in[O(G)]$ belongs to $O(G)$ if and
only if for each of the four possible pairs $(\e,\e')\in\{0,1\}\T\{0,1\}$,
both rational functions 
$$
f_{\jj^\e,\jj^{\e'};\pm}^*(\ph), f_{\jj^\e,\jj^{\e'};\mp}^*(\ph)\in[O(\CC^\nu\T T\T\CC^\nu)]
$$
belong to $O(\CC^\nu\T T\T\CC^\nu)$. 
\endproclaim

The proof of Theorem 5.1 will rely on the following statement.

\proclaim{Lemma 5.3} $\dim(G-((U^+TU^-)\cup(U^-TU^+)))\le\dim(G)-2$.
\endproclaim

\demo{Proof}
Using the Bruhat decomposition, we obtain: 
$$\align
G-((U^+TU^-)\cup(U^-TU^+))&=(G-(B^+U^-))\cap(G-(B^-U^+))\\
&=\Bigl(\bigcup_{w\in W-\{1\}}B^+\dw U^-\Bigr) \cap \Bigl(\bigcup_{w'\in W-\{1\}}B^-\dw'U^+\Bigr)\\
&=\bigcup_{w,w'\text{ in }W-\{1\}}(B^+\dw U^-)\cap(B^-\dw'U^+).\endalign$$
It is therefore enough to show that for any $w\ne 1$ and $w'\ne1$, we have
$$
\dim((B^+\dw U^-)\cap(B^-\dw'U^+))\le\dim(G)-2.
\tag 5.3.1
$$
This is clear if either $B^+\dw U^-$ or $B^-\dw'U^+$ has dimension $\le\dim(G)-2$. 
Thus we can assume that $\dim(B^+\dw U^-)=\dim(B^-\dw'U^+)=\dim(G)-1$ or equivalently $|w|=|w'|=1$.
Then both $\ov{B^+\dw U^-}$ and $\ov{B^-\dw'U^+}$
(closures in $G$) are irreducible of dimension $\dim(G)-1$.
If $\ov{B^+\dw U^-}\ne\ov{B^-\dw'U^+}$, then
$$\dim(\ov{B^+\dw U^-}\cap\ov{B^-\dw'U^+})\le\dim(G)-2,$$ 
implying (5.3.1).
Thus we may assume that $\ov{B^+\dw U^-}=\ov{B^-\dw'U^+}$.
By our assumption, $w=s_i$ for some $i\in I$.
For any $c\in\CC$ we have $y_i(c)B^-\dw'U^+\sub B^-\dw'U^+$ hence
$y_i(c)\ov{B^-\dw'U^+}\sub\ov{B^-\dw'U^+}$. Using our assumption, 
we also deduce that $y_i(c)\ov{B^+\ds_i U^-}\sub\ov{B^+\ds_iU^-}$ for any $c\in\CC$.
We have $B^+\ds_iU^-=B^+(s_iw_0)\dot{}\,U^+\dw_0\i$. 
For $c\in\CC^*$, we have
$$y_i(c)B^+\ds_iU^-\sub B^+\ds_iB^+B^+(s_iw_0)\dot{}U^+\dw_0\i
\sub B^+\dw_0B^+\dw_0\i=B^+U^-$$
and this is disjoint from $B^+\ds_iU^-$.
(We have used that $|s_i(s_iw_0)|=|s_i|+|s_iw_0|$.) This contradicts
the inclusion $y_i(c)\ov{B^+\ds_i U^-}\sub\ov{B^+\ds_iU^-}$.
\qed\enddemo

\demo{{\bf 5.4.} Proof of Theorem 5.2}
The ``only if'' statement in Theorem~5.2 is obvious. 
Let us prove the ``if'' statement. 
Consider $\ph\in[O(G)]$ such that the eight conditions in Theorem~5.2 are satisfied.
We need to show that $\ph\in O(G)$. 

Suppose that $G$ is of type $A_1$. Then
$\ph$ is regular on $U^+TU^- \cup U^-TU^+$. 
Hence by Lemma~5.3, it is regular on~$G$, and we are done.

In the rest of the proof, we assume that $G$ is of type other than~$A_1$.

We will first show that $\ph$ is a regular function on the open set
$U^+TU^-$. 

From our assumptions we see---as in the proof in 3.14---that
(using the same notation as~4.4)
$\ph$ is regular on each of the following open subsets of $U^+TU^-$: 

$\bullet$\ $\cv^\e_* T\cv^{\e'-}_*$, for $\e,\e'\in\{0,1\}$; 

$\bullet$\ $\cv^\e_* TU^-_{\jj^\e;i}$, for $\e\in\{0,1\}$ and $i\in I_{[\e+h+1]}$; 

$\bullet$\ $U^+_{\jj^\e;i}T\cv^{\e'-}_*$, for $\e,\e'\in\{0,1\}$ and $i\in I_{[\e+h+1]}$; 

$\bullet$\ $U^+_{\jj^\e;i}TU^-_{\jj^{\e'};i'}\,$, 
for $\e,\e'\in\{0,1\}$, $i\in I_{[\e+h+1]}$, and $i'\in I_{[\e'+h+1]}$.  

\noindent
Hence $\ph$ is regular on the union of these
subsets, i.e., on $\cu T{}\cu^-$ (where $\cu\sub U^+$ is given by (3.13.1)
and $\cu^-=\io(\cu)\sub U^-$). We have
$$U^+TU^- -\cu T\cu^- = ((U^+-\cu)T\cu^-) \cup (\cu T(U^- -\cu^-)).$$
By Lemma 3.13, we have
$$\dim((U^+-\cu)T\cu^-)\le\nu-2+r+\nu=\dim(G)-2
$$
and similarly 
$$\dim(\cu T(U^- -\cu^-))\le\dim(G)-2.
$$
It follows that
$$\dim(U^+TU^- -\cu T\cu^-)\le\dim(G)-2.
$$
Since $\ph$ is regular on $\cu T\cu^-$, we conclude that $\ph$ is regular
on $U^+TU^-$.
An entirely similar argument shows that $\ph$ is regular on $U^-TU^+$.
It follows that $\ph$ is regular on the open subset
$(U^+TU^-)\cup(U^-TU^+)$ of $G$. 
Together with Lemma~5.3,
this implies that $\ph$ is regular on $G$. 
\qed\enddemo


\widestnumber\key{ABCDE}
\Refs
\ref\key{BFZ96}\by A. Berenstein, S. Fomin and A. Zelevinsky
\paper Parametrizations of canonical bases and totally positive matrices
\jour Adv. Math. \vol 122 \yr 1996 \pages 49--149\endref
\ref\key{BZ97}\by A. Berenstein and A. Zelevinsky\paper Total positivity in Schubert varieties
\jour Comm. Math. Helv.\vol72\yr1997\pages128--166\endref
\ref\key{Bo}\by N. Bourbaki  \book Groupes et alg\`ebres de Lie, Chapitres 4--6
\publ Hermann\publaddr Paris
\yr1968\endref
\ref\key{FZ99}\by S. Fomin and A. Zelevinsky\paper Double Bruhat cells and total positivity
\jour Jour. Amer. Math. Soc.\vol12\yr1999\pages335--380\endref
\ref\key{KL79}\by D. Kazhdan and G. Lusztig\paper Representations of Coxeter groups and Hecke algebras\jour
Inv. Math.\vol53\yr1979\pages 165--184\endref
\ref\key{L94}\by G. Lusztig\paper Total positivity in reductive groups\inbook Lie theory and geometry\bookinfo
Progr. in Math.\vol123 \publ Birkh\"auser \yr1994\pages 531--568\endref
\ref\key{L98}\by G. Lusztig\paper Introduction to total positivity\inbook Positivity in Lie theory: open
problems\bookinfo  ed. J. Hilgert et al.\publ de Gruyter\yr1998\pages 133--145\endref
\ref\key{L19}\by G. Lusztig\paper Total positivity in reductive groups,II\jour Bull. Inst. Math. Acad. Sinica
\vol14\yr2019\pages403--460\endref 
\endRefs
\enddocument